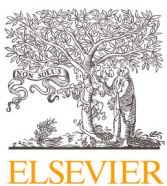
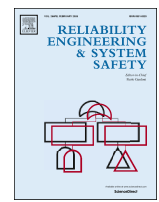

# Fortifying critical infrastructure networks with multicriteria portfolio decision analysis: An application to railway stations in Finland

Joaquín de la Barra [*], Ahti Salo, Leevi Olander, Kash Barker, Jussi Kangaspunta

ABSTRACT

Advanced societies are crucially dependent on critical infrastructure networks for the reliable delivery of essential goods and services. Hence, well-founded analyses concerning disruptions are necessary to inform decisions that aim to ensure the performance of these networks in the face of failures caused by vulnerabilities to external hazards or technical malfunctions. In this setting, we develop an approach based on multicriteria decision analysis to support the identification of cost-efficient portfolios of preventive fortification actions. Our approach (i) accounts for multiple performance objectives, such as those that maximize the uninterrupted volume of traffic between different origin-destination pairs in a transportation network, (ii) uses methods of probabilistic risk assessment to quantify the expected performance of the network with regard to these objectives, and (iii) uses a search algorithm combined with an optimization model to identify those combinations of fortification actions that are cost-efficient in improving the performance of the network, given the available, possibly partial information about the relative importance of objectives and minimum performance requirements on them. Our methodological contributions are illustrated by a case study on the analysis of railway switches at a representative Finnish railway station.

## 1. Introduction

Critical infrastructures comprise all the assets, systems, and networks that provide essential functions for society [1]. These infrastructures are essential for ensuring the continuity of operations in sectors such as energy, water supply, transportation, and telecommunications. Because disruptions in these infrastructures can cause significant damage to public health, safety, security, and the economy, they must function adequately to enable the attainment of goals related to economic productivity, sustainability, and social development [2].

In Europe, railway networks amounted to 202,000 km in 2022, with considerable recent growth in high-speed trains [3]. Disruptions in railway networks can undermine essential performance objectives, such as ensuring connectivity between strategic locations or providing reliable connections for the transport of passengers and goods. These disruptions can be caused by failures resulting from the normal wear and tear of technical systems or by vulnerabilities to external hazards, including extreme weather conditions and intentional attacks. Consequently, there is a need to understand what types of disruptions can cause damage to the railway network, how these disruptions can affect transportation volumes, and what fortification actions are cost-efficient in mitigating the harmful impacts as measured by the reduced volume of transportation, subject to the constraint that there are limited resources [see, e.g. 4–6].

Not all components in transportation networks are equally important, because disruptions in some parts of the network erode performance more than those in other parts. Moreover, because impacts on network performance depend on the states of all network components [see, e.g. 7], there is a need to evaluate the impacts of *combinations* of disruptions (e.g., a single component may not be important on its own; but if it fails simultaneously with some other component, there may be drastic consequences for network performance). Furthermore, because disruptions are inherently uncertain, one needs to analyze how the expected impacts of disruptions depend on alternative assumptions about the underlying probabilities of component failures, also in situations where there is partial information about these probabilities [8].

Fortification actions lower the probability of disruptions that typically undermine the performance of the network. Because there are usually many fortification actions that are implemented together, we analyze *portfolios* which consist of many actions. Our analyses support the cost-efficient allocation of resources to those portfolios of fortification actions that most improve the expected performance relative to the

---

costs of implementing these actions [9]. For example, these analyses can be produced for different budgets to prepare plans for the cost-efficient allocation of additional resources, or, in the event of budget cuts, to identify which previously planned fortification actions can be canceled with minimal reduction in expected network performance.

In this paper, we develop an approach based on multicriteria portfolio decision analysis in order (i) to assess and aggregate multiple objectives that reflect the services provided by infrastructure networks and (ii) to guide the cost-efficient allocation of resources to fortification actions that are focused on the components of the network, with the aim of providing the maximum expected performance level of its services. Technically, we represent the network as a graph consisting of nodes (components that may fail) and edges (fully reliable connections between components), whereby our formulation helps identify which nodes, either individually or jointly, are most important to the network's performance. Importantly, this formulation identifies which portfolios of fortification actions are cost-efficient in enabling the targeted level of network performance.

The remainder of this paper is organized as follows. Section 2 discusses earlier approaches to analyzing disruptions and their impacts on the performance of transportation networks. Section 3 develops our approach, which combines probabilistic risk assessment with multicriteria portfolio decision analysis to quantify the expected performance of the network with respect to several objectives and, in addition, guides the allocation of resources to those actions that fortify the network cost-efficiently. Section 4 presents a case study on the fortification of railway switches to improve the reliability of connections at a representative railway station in Finland. Section 5 discusses the numerical results, and Section 6 outlines avenues for future research.

## 2. Background

Early studies of railway networks focused mainly on topological configuration due to limited data and computational power [10]. Later, there has been a proliferation of specific models to reduce travel times [11], design schedules [12], plan new lines [13], reduce operating costs [14], and ensure track functionality, among others. In many countries, infrastructure decisions (e.g., ensuring the reliability of the railway network) and operational decisions (e.g., developing train schedules) are made by different entities. In part due to this separation of responsibilities, research activities supporting decisions on the development and maintenance of railway infrastructure have historically been somewhat decoupled from efforts to inform decisions in the operational planning and management of transportation services. In this paper, we adopt a general view of the network as a critical infrastructure that enables essential services.

Many analyses of railway systems employ network theory to identify critical components, evaluate network performance, and develop strategies to reinforce or disrupt them [see, e.g., 15–18]. Latora and Marchiori [19] present a method to find the critical components of an infrastructure network represented by nodes and edges. They also analyze how improvements, such as adding edges between nodes, enhance network performance. In their model, performance is measured using a topological metric, which utilizes only information about the nodes' positions and connections in the network to quantify the efficiency of information exchange over the network. Ip and Wang [20] also propose a methodology to evaluate the resilience of railway networks based on topological metrics, such as the number of independent paths. They also propose an algorithm to fortify the network to attain improved resilience. Fecarotti et al. [21] propose a non-linear integer programming model that considers topological metrics to select maintenance strategies.

Although several studies have used topological metrics to assess network performance, few have evaluated the quality of these assessments [22]. In this regard, a general drawback of topological metrics is that they may not be closely related to network performance objectives. For example, the multi-objective optimization approach proposed by Hao et al. [23] helps identify critical components in a network, taking into account its interactions with other systems. Yet, one of their key findings is that the criticality of the nodes is not necessarily related to topological importance. Based on a comparison of a topological metric of network nodes with an importance metric derived from the objective reflecting the performance of the network, Olander et al. [24] found only a weak correlation between the two, suggesting that the topological metric may not be suitable for estimating the importance of the nodes. Based on their analysis of the robustness of power systems in different disruption scenarios, LaRocca et al. [25] conclude that many topological metrics are of limited value. The criticality of a node is a system-level phenomenon in that the criticality of a given node can depend more on what other nodes are disrupted rather than on the component itself, as shown by Alderson et al. [26]

A further challenge in evaluating the performance of networks under uncertainty is identifying possible hazards and their impacts on the network and its components. In some cases, it is unclear what these hazards are [27] or to what extent the infrastructure will withstand expected or unexpected hazards [28]. Zhang et al. [29] summarize recent studies on quantifying the loss of functionality of railway systems due to various hazards, such as extreme weather or seismic events. They propose a framework to simulate the loss of functionality due to combined hazards in a coupled railway and airline system, noting that modeling the propagation of disruptions is an essential step in quantifying this loss. It is worth noting that the assumption about independent failure events can lead to overestimating network reliability, resulting in poor design and unacceptable performance. Nazarizadeh et al. [30] propose a model that includes common cause failures and interactive failures. They present a case study on the Iranian railway system, which demonstrates that their model can provide improved reliability estimates.

Several authors describe scenarios associated with the uncertain realization of hazards and their respective probability estimates. For example, Joshi et al. [31] and Yang et al. [32] consider scenarios of rainfall and tornadoes to assess the risks affecting the railway systems in India and China, respectively. Turoff et al. [33] propose a collaborative dynamic scenario model based on expert judgments to estimate the cascading effects of critical infrastructure interactions.

Zio [28] and Sedghi et al. [34] call for the development of frameworks to help railway infrastructure managers understand and quantify the complexity of the network and, by doing so, help them prepare for hazards to ensure acceptable network performance.

In this context, we develop a general assessment approach and accompanying search algorithms to assess the importance of nodes in a network, considering the possibility of node disruptions, as reflected by node-specific disruption probabilities; multiple objectives associated with the different services that are enabled by the network; and possible minimum performance requirements on the attainment of these objectives. These assessments of node importance are further leveraged to offer recommendations for fortification actions that help improve the expected performance level of the network. The referenced models and our approach are summarized in Table 1. The works of [18] and [15] are not included, as they only provide literature reviews on Bayesian networks and graph theory, respectively.

## 3. Proposed approach for fortifying networks

Fig. 1 illustrates the overall structure of the proposed approach. Information about the services enabled by infrastructure is used to derive corresponding weights that reflect the importance of these services. Information about physical network infrastructure is used to define the network topology and assess the impacts of node disruptions on the ability to enable requisite services. Information about decision alternatives consists of the set of available fortification actions, characterized by the cost of implementation and the impact on lowering the probability of failure at nodes. These three sources of information are integrated within an optimization framework to identify cost-efficient portfolios of fortification actions. The following sections provide a detailed description of the model.

**Table 1**
Summary of related approaches. SA: search algorithm, INLP: nonlinear integer programming, MOO: multi-objective optimization, PRA: probabilistic risk assessment.

| Ref. | Goal | Methodology |
| --- | --- | --- |
| [20] | Maximize performance | Topological metrics, SA, MOO |
| [35] | Maximize performance, identify important components | Topological metrics |
| [21] | Maximize performance | Topological metrics, INLP, PRA |
| [24] | Identify important components | Topological metrics, performance-based metrics |
| [17] | Identify important components | Topological metrics |
| [26] | Identify important components | Optimization |
| [23] | Identify important components | Topological metrics, SA, MOO |
| [25] | Risk assessment | Topological metrics, physical flow models |
| [29] | Risk assessment | Simulation, optimization |
| [16] | Risk assessment | Spanning trees |
| [30] | Risk assessment | Reliability block diagrams, dependent failures |
| [33] | Risk assessment | Failure scenarios, Delphi method |
| [32] | Risk assessment | Failure scenarios, C-vine copula |
| [31] | Risk assessment | Failure scenarios |
| This work | Maximize performance, identify important components | SA, PRA, Information sets |

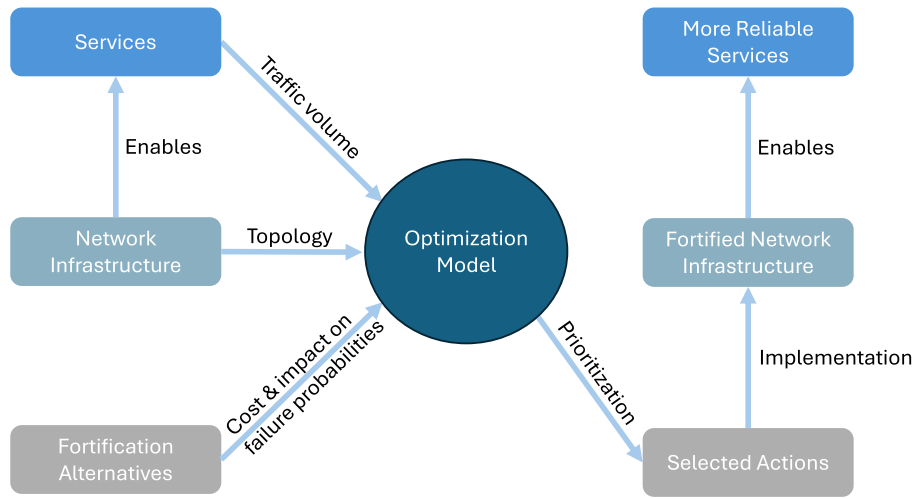

**Fig. 1.** Overall structure of the proposed approach.

### 3.1. Network representation

Let $G(V; E)$ denote a network consisting of nodes $V = \{1, \ldots, n\}$ and undirected edges $E \subseteq \{(i, i') \mid i, i' \in V\}$ between the nodes. In our analysis, a path is an acyclic sequence of edges that connect two nodes in the network. The state of a node is binary; either operational or disrupted. If a node is disrupted, none of the paths containing it can be traversed. Thus, if the nodes $D \subseteq V$ are disrupted, the remaining network is $G(V^O; E^O)$ such that $V^O = V \setminus D$ and $E^O = \{(i, i') \in E \mid i, i' \in V^O\} \subseteq E$.

Uncertainty about the state $x_k$ of node $k$ is modeled as a binary random variable $X_k$ such that $x_k = 0$ if node $k \in V$ is disrupted and $x_k = 1$ if it is operational. The state of the network $x = (x_1, \ldots, x_n) \in \mathcal{X} = \{0, 1\}^n$ is a combination of states for all nodes. Thus, there are $2^n$ network states. We assume that node disruptions occur independently of each other, i.e., the failure probability of a node does not depend on the states of other nodes. The probability of disruption at node $k$ is denoted by $p_k = \mathbb{P}[X_k = 0]$. Due to the independence assumption, the probability distribution over all network states is implied by the vector $p = (p_1, \ldots, p_n)$. The probability that the network is in state $x \in \mathcal{X}$ is denoted by $p(x)$.

### 3.2. Assessing network performance

Because critical infrastructure networks often enable multiple services, their performance can be measured with several objectives. For example, the performance of a transportation network can be assessed in terms of the reliability of connections between relevant origin-destination pairs. In the literature, the objectives have also been assessed with topological metrics, such as the average of the shortest distances between network nodes [see, e.g., 19,36].

The performance of the network depends on its state, as node disruptions usually lower the reliability of services. Concretely, in the case of $m$ services, there is an equal number of objectives whose attainment can be measured by introducing corresponding criteria and, more specifically, by utility functions $u_j(\cdot) : \mathcal{X} \mapsto [0, 1]$ such that $u_j(x)$ represents the utility that is associated with the performance on criterion $j = 1, \ldots, m$ when the network is in state $x \in \mathcal{X}$.

We assume that the network is coherent, so that these utilities do not increase if a node fails. The utilities are normalized so that the lowest utility $u_j(x) = 0$ on criterion $j$ is attained when all nodes are disrupted $x_k^\circ = 0, k = 1, \ldots, n$. Conversely, the highest utility $u_j(x) = 1$ is attained when all nodes are operational $x_k^* = 1, k = 1, \ldots, n$.

The utility functions $u_j(\cdot)$ can be aggregated with an additive multi-criteria utility function in Eq. (1) if criteria are mutually preferentially independent and every criterion is additive independent (for details, see [37]). Specifically, in the utility function

$$u(x, w) = \sum_{j=1}^{m} w_j u_j(x) \in [0, 1], \quad (1)$$

the weight of the criterion $w_j \in [0, 1]$, $j = 1, \ldots, m$ can be interpreted as the increase in the overall utility that can be attributed to criterion $j$ when the network state improves from its worst state (all nodes are disrupted) to the best one (all nodes are operational). Following the usual convention, these weights can be normalized so that $\sum_{j=1}^{m} w_j = 1$.

However, it can be challenging to elicit complete preference information in terms of exact criteria weights that represent the relative importance of the criteria. For example, if the criteria represent the planned volume of traffic between different origin-destination pairs in the network, this volume may not be fully known at the time of choosing fortification actions. Furthermore, decision-makers (DM) may be hesitant to make judgments about the importance of criteria, or there can be several DMs with different priorities.

In response to this recognition, and in order to support extensive sensitivity analyses, we characterize information about the relative importance of criteria with an *information set* $S$ which consists of all the weights that are compatible with the elicited preference statements [38]. These statements can be elicited by asking the DM to express their views on the relative importance of the criteria. Thus, for example, if criterion 1 is at least as important as criterion 2 but not more than two times more important, the constraints $w_2 \leq w_1 \leq 2w_2$ hold. As most types of such preference statements correspond to linear constraints, the resulting weight set $S$ – which is a subset of the non-informative weight set $S^0$ – is given by

$$S = \{w \in \mathbb{R}^m \mid Aw \leq b\} \subseteq \left\{w \in \mathbb{R}^m \mid w_j \geq 0 \, \forall j, \sum_{j=1}^{m} w_j = 1\right\} = S^0, \quad (2)$$

where the constraint matrix $A \in \mathbb{R}$ and the vector $b \in \mathbb{R}$ contain the coefficients of the constraints.

### 3.3. Fortifying the network

The reliability of the network can be improved through fortification actions that lower the probabilities of node disruptions. Fig. 2 illustrates the decision to implement an action to fortify a single network node. Without this action, the probability of disruption of node $k$ is $p_k$. If the action is implemented, this probability is reduced to $p'_k < p_k$.

Typically, there are many possible fortification actions that can be implemented to improve the reliability of the network. As a result, the network reliability depends on the *portfolio* (i.e., collection) of actions that are implemented. We assume that there are $h$ possible fortification actions, represented by the binary variables $q_1, \ldots, q_h$ such that $q_l = 1$ if action $l$ is implemented and $q_l = 0$ if not, $l = 1, \ldots, h$. Thus, the selected portfolio of fortification actions is a vector $q = (q_1, \ldots, q_h) \in \mathcal{Q} = \{0, 1\}^h$. In what follows, the total cost $c_{tot}(q) = \sum_{l=1}^{h} c_l q_l$ of portfolio $q$ is the sum of the costs of the individual actions $c_l$ that are contained in it. More complex cost structures can be introduced to represent cost synergies between actions, for instance.

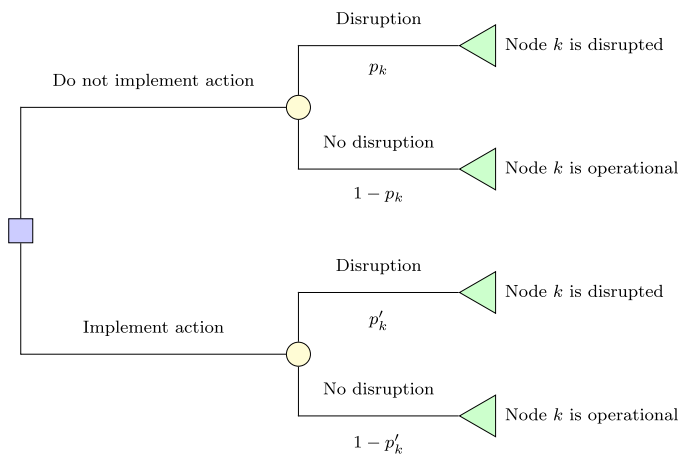

**Fig. 2.** Decision tree for implementing an action to fortify network node $k$.

A portfolio of fortification actions $q$ is *feasible* if (i) its total cost $c_{tot}(q)$ does not exceed the budget level of available resources $r$ (i.e., $c_{tot}(q) \leq r$) and (ii) satisfies relevant logical constraints (e.g., if actions 1 and 2 are mutually exclusive, then only one of them can be implemented, and thus constraint $q_1 + q_2 \leq 1$ must hold). The set of feasible portfolios is denoted by $\mathcal{Q}_F \subseteq \mathcal{Q}$. The aim is to determine all feasible portfolios that satisfy the relevant constraints and contribute most to the attainment of objectives, such as maximizing the reliability of the services enabled by the network.

The probability of disruption $p_k(q)$ at node $k$ depends on the portfolio $q \in \mathcal{Q}_F$ of actions. Thus, assuming that the probabilities of node disruptions depend on these actions only, the network state $x \in \mathcal{X}$ occurs with probability

$$p(x \mid q) = \prod_{k=1}^{n} \left[x_k(1 - p_k(q)) + (1 - x_k)p_k(q)\right] \in [0, 1]. \quad (3)$$

This formulation is general in that the probability of any node disruption can depend on the entire portfolio of implemented actions. Thus, if two actions lower the probability of disruption at a given node and these two actions can be taken jointly, one would need to estimate the probability of disruption at this node for situations in which only one of the alternative actions or both actions are implemented.

Because we are analyzing the reliability of services enabled by operational connections in a transportation network, the impact of fortification actions on the probabilities of providing operational connections must first be established. The next step is to identify which portfolios most improve the reliability of these connections. When the weights $w$ that represent the importance of connections are specified exactly, for instance, based on transportation volumes, the optimal portfolio can be obtained by solving the optimization model in Appendix A, to maximize the expected value of the utility in Eq. (1). On the other hand, when these weights are characterized by information sets that contain many weights, analytical concepts such as non-dominance and cost-efficiency can be used to compare portfolios of fortification actions.

### 3.4. Non-dominated and cost-efficient portfolios

In the case of information sets that contain many feasible weights, the aim is to determine which feasible portfolios outperform others in view of all these feasible weights. This can be done by comparing portfolios based on the concept of *dominance*.

**Definition 1.** Portfolio $q^1 \in \mathcal{Q}_F$ is dominated by portfolio $q^2 \in \mathcal{Q}_F$ in the information set $S$, denoted by $q^2 \overset{S}{\succ} q^1$, if and only if $\mathbb{E}\left[u(x, w) \mid q^1\right] \leq \mathbb{E}\left[u(x, w) \mid q^2\right]$ for all $w \in S$ and (ii) $\mathbb{E}\left[u(x, w) \mid q^1\right] < \mathbb{E}\left[u(x, w) \mid q^2\right]$ for some $w \in S$.

The dominance between two portfolios can be determined by comparing the expected utility of these portfolios at the extreme points $w_e \in S^{ext}$ of the information set $S$ [39]. This information set is a polyhedral set, and thus these points can be computed with linear programming techniques [see, e.g., 40]. Moreover, if $\mathbb{E}\left[u(x, w) \mid q^1\right] = \mathbb{E}\left[u(x, w) \mid q^2\right]$ $\forall w \in S$, portfolios $q^1$ and $q^2$ have the same expected utility, which is denoted by $q^1 \overset{S}{\sim} q^2$.

A feasible portfolio is cost-efficient if (i) it is not dominated by any feasible portfolio of equal or lower cost and (ii) there is no other portfolio with lower cost with equal expected utility.

**Definition 2.** Portfolio $q^1 \in Q_F$ is cost-efficient with respect to another portfolio $q^2 \in Q_F$ in the information set $S$, denoted by $q^1 \succ_C^S q^2$, if and only if (i) $q^1 \succ^S q^2, c(q^1) \leq c(q^2)$ or (ii) $q^1 \sim^S q^2$ and $c(q^1) < c(q^2)$.

**Definition 3.** Portfolio $q^1 \in Q_F$ is cost-efficient in the information set $S$, denoted by $q^1 \in Q_{CE} \subseteq Q_F$ if and only if $\nexists q^2 \in Q_F$ such that $q^2 \succ_C^S q^1$.

Although Definitions 1–3 refer to the maximization of expected utility, these definitions can be extended to account for feasibility constraints that may be derived from risk measures such as value at risk (VaR) and conditional value at risk (CVaR) [41].

### 3.5. Determination of non-dominated portfolios

#### 3.5.1. Algorithm to determine cost-efficient portfolios

Algorithm 1 determines the cost-efficient portfolios $Q_{CE}$. At each iteration, it generates new portfolios by adding a fortification action to the previously computed set of cost-efficient portfolios and the *basic portfolios* (i.e., portfolios that are not cost-efficient but can potentially be extended to cost-efficient portfolios by adding fortification actions). In particular, portfolios that cannot be in the set of basic portfolios are removed to avoid enumerating all feasible portfolios.

---

**Algorithm 1** Compute $Q_{CE}$.

1: Compute $u(x, w)$ for all $x \in \mathcal{X}$ and $w \in S^{ext}$
2: $Q^0 \leftarrow \{(0, \ldots, 0)\}$
3: $Q_B^0 \leftarrow \{\}$
4: **for** $l = 1$ to $h$ **do**
5: $\quad Q^l \leftarrow \left\{ q^1 \in Q_F \mid q_l^1 = 1 \wedge \exists q^2 \in Q_B^{l-1} \cup Q^{l-1} : q_y^1 = q_y^2, \forall y \neq l \right\}$
6: $\quad Q_D^1 \leftarrow \left\{ q^1 \in Q^l \mid \exists q^2 \in Q^l \cup Q^{l-1} : q^2 \succ_C^{S^{ext}} q^1 \right\}$
7: $\quad Q^l \leftarrow Q^l \setminus Q_D^1$
8: $\quad Q_D^2 \leftarrow \left\{ q^1 \in Q^{l-1} \mid \exists q^2 \in Q^l : q^2 \succ_C^{S^{ext}} q^1 \right\}$
9: $\quad Q^{l-1} \leftarrow Q^{l-1} \setminus Q_D^2$
10: $\quad Q^l \leftarrow Q^l \cup Q^{l-1}$
11: $\quad Q_B^l \leftarrow \left\{ q \in Q_D^1 \cup Q_D^2 \cup Q_B^{l-1} \mid \nexists q^1 \in Q^l : q^1 \succ_C^{S^{ext}} q^a \right\}$
$\quad\quad$ with $q_y^a = q_y, \forall y \leq l \wedge q_y^a = 1, \forall y > l$
12: **end for**
13: $Q_{CE} \leftarrow Q^h$

---

In Step 1, the utility function is computed for all network states and all extreme points of the weight set. In Step 2, the set of cost-efficient portfolios is initialized by including only the empty portfolio, which has no fortification actions. In Step 3, the set of *basic portfolios* is initially empty. In Steps 4 to 12, the index $l = 1, \ldots, h$ iterates through all fortification actions. In Step 5, a set $Q^l$ is built by taking every portfolio in the set $Q^{l-1} \cup Q_B^{l-1}$ and modifying its $l$th action. In Step 6, the portfolios in $Q^l$ are compared to portfolios in $Q^{l-1}$, and cost-inefficient portfolios are stored in the set $Q_D^1$. In Step 7, the cost-inefficient portfolios in $Q_D^1$ are removed from the set $Q^l$. In Step 8, the portfolios in $Q^{l-1}$ are compared to portfolios in $Q^l$, and cost-inefficient portfolios are stored in the set $Q_D^2$. In Step 9, the cost-inefficient portfolios in $Q_D^2$ are removed from the set $Q^{l-1}$. In Step 10, the set $Q^l$ is updated, including the cost-efficient portfolios from $Q^{l-1}$. In Step 11, the set of basic portfolios is updated by retaining the previous basic portfolios and introducing those cost-inefficient portfolios that can be extended to cost-efficient portfolios. To determine if a portfolio $q$ is introduced in the set of basic portfolios, the algorithm constructs its extended portfolio $q^a$ by adding all the remaining fortification actions to $q$. If the extended portfolio $q^a$ of a portfolio $q$ is not cost-efficient, then $q$ is not in the set of basic portfolios. The analysis of the extended portfolio $q^a$ provides an upper bound of the expected utility that can be achieved by extending a portfolio. Still, they are not necessarily feasible (e.g., their cost may exceed the budget). If a portfolio $q^a$ is infeasible, the upper bound does not accurately reflect the value of the maximum expected utility. As a result, the algorithm may retain in the set basic portfolios some portfolios that need to be discarded later on. In B, we provide an alternative approach to generate tighter bounds to remove further portfolios. The algorithm terminates in Step 13, returning the set of cost-efficient portfolios $Q^h$.

#### 3.5.2. Binary utility function for individual objectives

The utility functions in (4) have binary values to represent whether or not the network fulfills a given condition. For example, in a transportation network, the objective may refer to the existence of operational connections for evacuation purposes or the transportation of sensitive cargo. In such cases, the resulting objective can be represented by the binary utility function

$$u_j(x) = \begin{cases} 1, & \text{if objective } j \text{ is met in network state } x \in \mathcal{X}. \\ 0, & \text{otherwise.} \end{cases} \quad (4)$$

In this binary case, the expected utility with regard to this single-objective utility function is the probability with which this objective is attained, which is also the reliability of this objective, i.e.,

$$\mathbb{E}[u_j(x)] = 1 \cdot \mathbb{P}[u_j(x) = 1] + 0 \cdot \mathbb{P}[u_j(x) = 0] = \mathbb{P}[u_j(x) = 1]. \quad (5)$$

#### 3.5.3. Minimum performance requirements for individual objectives

In some cases, one may require that the reliability of objective $j$ is greater than or equal to $\alpha_j$. This requirement can be represented by the chance constraint $\mathbb{E}[u_j(x)] \geq \alpha_j$. In the case of many objectives, the vector of those binary performance requirements is denoted by $\alpha = (\alpha_1, \ldots, \alpha_m)$, where $\alpha_j = 0$ if there is no such requirement.

Algorithm 1 does not consider such minimum performance requirements for individual objectives because it does not consider the utility functions for individual objectives directly. Instead, it uses the overall utility function of the network performance in Eq. (1), which is a weighted sum of utility functions for the objectives.

The non-informative set $S^0$ in Eq. (2) is a special case as it imposes no constraints on the importance of individual objectives. The extreme points of the set $S^0$ are the canonical vectors $e_j = (w_1, \ldots, w_j, \ldots, w_m) \in S^{0, ext}$, where $w_j = 1$ and $w_i = 0$ for all $i \neq j$. At every extreme point $e_j \in S^{0, ext}$, the utility function of the network is the single-objective utility function, i.e., $u(x, e_j) = u_j(x)$. Thus, for each objective, there is a cost-efficient portfolio in the non-informative set that gives the maximum performance on this objective. This is formalized in Lemma 1.

**Lemma 1.** $\forall j = 1, \ldots, m \; \exists q^1 \in Q_{CE}$ in the information set $S^0$ such that $\mathbb{E}[u_j(x) \mid q^1] \geq \mathbb{E}[u_j(x) \mid q^2] \; \forall q^2 \in Q_F \setminus Q_{CE}$.

**Proof.** We prove Lemma 1 by contradiction. Consider the objective $j$. Let $q^1 \in Q_{CE}$ denote the portfolio with the maximum expected performance in $j$. Assume, for the sake of contradiction, that there exists $q^2 \in Q_F \setminus Q_{CE}$ such that $\mathbb{E}[u_j(x) \mid q^2] > \mathbb{E}[u_j(x) \mid q^1]$. Since $q^1$ is chosen as the portfolio with the maximum expected performance in $j$, it follows that $\nexists q \in Q_{CE}$ such that $\mathbb{E}[u(x, e_j) \mid q^2] < \mathbb{E}[u(x, e_j) \mid q]$, where $e_j \in S^0$ is a canonical vector.

By Definitions 1 and 2, this implies that $q^2$ is non-dominated in $S^0$. Consequently, $q^2 \notin Q_F \setminus Q_{CE}$, contradicting the initial assumption. □

The maximum performance on objective $j$ is attained for some cost-efficient portfolio in the information set $S$ that contains the vector $e_j \in S$. Therefore, Algorithm 1 for determining the cost-efficient portfolios can be adapted in the case of minimum performance requirements for individual objectives, provided that the corresponding extreme point is in the information set ($\alpha_j \neq 0 \Rightarrow e_j \in S^{ext}, \; \forall j = 1, \ldots, m$). Once the cost-efficient portfolios have been determined, the set must be filtered to remove those that do not fulfill requirements on expected performance.

#### 3.5.4. Algorithm to determine cost-efficient portfolios with minimum performance requirements

Algorithm 2 computes cost-efficient portfolios that fulfill minimum performance requirements. To achieve this, the information set is expanded to encompass all extreme points that meet the minimum performance requirements.

---

**Algorithm 2** Compute $Q_{CE}$ with minimum performance requirements $\alpha$.

1: $S^\alpha \leftarrow \{e_j \mid \alpha_j \neq 0 \ \forall j = 1, \ldots, m\}$
2: $S^* \leftarrow S^{ext} \bigcup S^\alpha$
3: $Q^*_{CE} \leftarrow$ Algorithm 1 with $S^*$
4: $Q^\alpha_{CE} \leftarrow \{q \in Q^*_{CE} \mid \mathbb{E}[u_j(x) \mid q] \geq \alpha_j \ \forall j = 1, \ldots, m\}$
5: $Q_{CE} \leftarrow Q^\alpha_{CE} \setminus \left\{ q^1 \in Q^\alpha_{CE} \mid \exists q^2 \in Q^\alpha_{CE} : q^2 \succ^{S^{ext}}_C q^1 \right\}$

---

In Step 1, the weight set $S^\alpha$ contains the canonical vector $e_j$ for every minimum performance requirement $\alpha_j$. In Step 2, the extended set $S^*$ is obtained by combining $S^\alpha$ with the extreme points of the set that represents the DM's preferences $S^{ext}$. In Step 3, the cost-efficient portfolios in the extended information set $S^*$ are determined with Algorithm 1. In Step 4, portfolios that do not fulfill the minimum performance requirements are removed. In Step 5, the cost-efficient portfolios in the information set $S^{ext}$ are computed. The set of portfolios $Q^\alpha_{CE}$ is usually much smaller than the set of feasible portfolios, so it is possible to use algorithms such as the one proposed by [42] to identify the set $Q_{CE}$.

## 4. Case study

We illustrate our approach by considering the Siilinjärvi railway station in Northern Savonia, Finland, depicted in Fig. 3(a). A network representation based on this station structure is in Fig. 3(b). The small squares represent the 22 railway switches, mechanical devices that allow trains to move from one railway track to another. Specifically, these switches correspond to the network's nodes, while the railway tracks connecting them correspond to its edges. Nodes A, B, and C are the border nodes that define the geographic boundaries of the station.

The station serves the commuting needs of the municipality of Siilinjärvi, as well as significant passenger and freight traffic from the south, east, and west. A precondition for enabling this traffic is that trains can pass through the station. The trains change railway tracks to go in different directions. For example, a train from the south can enter the station at border node C, exit at border node B, and then proceed east. This is possible only if there are operational switches that provide at least one operational path from C to B. Due to the network structure, such a path may exist even if some switches are disrupted.

The railway tracks between switches are far less prone to disruptions than switches, which are relatively complex multi-component systems that may fail, for example, due to technical failures or adverse weather conditions [43]. From this perspective, it is of primary interest to determine which switches are key for ensuring the reliability of connections between the three pairs of directions defined by the three border nodes.

The probabilities of switch failures depend on several factors, including age and technology. In practice, the estimation of these probabilities can be facilitated by condition-based monitoring systems, for example [44]. In this example, we assume that the switches are equally prone to failure and, specifically, that each switch fails with $p_k = 0.01$ during the planned maintenance period. Furthermore, possible fortification actions are equally costly as well as equally effective in lowering the probability of switch failure by 50 % to $p'_k = 0.005$. These simplifying assumptions, made only for the purposes of illustration, make it easier to show how the preference information, derived from the planned volume of traffic for the three different connections, is reflected in the recommended allocation of resources to those fortification actions that improve the reliability of especially those connections with a high traffic volume (see Fig. 5).

Table 2 shows how the physical parts of the railway station correspond to elements of the graph and the optimization model. The pro-posed portfolio optimization approach can be applied to other infrastructures with an analogous correspondence.

**Table 2**
Relation between real elements and their modeling.

| Model Element | Real Element/Explanation |
|---|---|
| Network $G(V; E)$ | Railway station |
| Node $i \in V$ | Railway switch |
| Edge $(i, i') \in E$ | Railway track connecting two switches |
| State of a node $x_i$ | Operational status of a switch (i.e., operational or disrupted) |
| Network state $x = (x_1, \ldots, x_n)$ | Vector of states for all switches, defining the overall system state |
| Fortification action $q_i$ | Maintenance action which lowers the failure probability at switch $i$ |
| Portfolio $q = (q_1, \ldots, q_h)$ | Set of fortification actions |

### 4.1. Assumptions

To summarize, this numerical case study is based on the following assumptions. First, uncertainties about node failures are characterized as probabilities that could be derived, for example, from condition-based monitoring systems. Second, failures occur only at the nodes, as edges are assumed to be fully reliable (a plausible assumption in our case, where edges represent railway tracks connecting switches). Third, node failures occur independently of each other. Fourth, each node has two possible states: operational or disrupted. Fifth, although the methodology allows for any probabilities of node failure, we consider equal probabilities and also equally costly fortification actions that are equally effective in lowering the probability of node failures. Sixth, the cost of a portfolio of fortification actions is given by the sum of the costs of the actions contained in it. Seven, although the methodology allows for considering multiple fortification actions per node, we consider a single fortification action. These last three assumptions are made only for the sake of illustration to demonstrate how preference information derived from the planned traffic volume influences recommendations for choosing fortification actions.

### 4.2. Station performance

We evaluate the station based on the reliability of the connections between pairs of border nodes A, B, and C. These connections are bidirectional, named as (A, B), (A, C), and (B, C). A connection (X,Y) is operational if there is at least one operational path between X and Y such that all switches on this path are operational. The reliability of a connection is the probability that it is operational. For each connection $j$, there is an associated binary utility function (4) which represents the objective "connection $j$ is operational." These three utility functions are included in the additive utility function (1) that measures the overall performance of the station.

The task of evaluating the utility function (4), which represents the objective of having an operational connection between a pair of border nodes, is known in the literature as the *terminal-pair reliability problem* [45]. This problem is NP-hard and, therefore, cannot be solved in polynomial time [46]. For small networks, it can be solved by determining if a path exists between the pair of border nodes for each network state. When the number of nodes $n$ is large, this is impractical because there are $2^n$ network states.

To reduce the computational complexity of the terminal-pair reliability, we employ the minimum cut upper bound approximation [47], acknowledging that other algorithms exist for computing the exact reliability even in large networks (see [48] for a survey of recent exact algorithms). If the network size cannot be handled by exact algorithms, or the computation time is critical, one can use approaches such as decomposition approaches [49], Monte Carlo simulation in combination

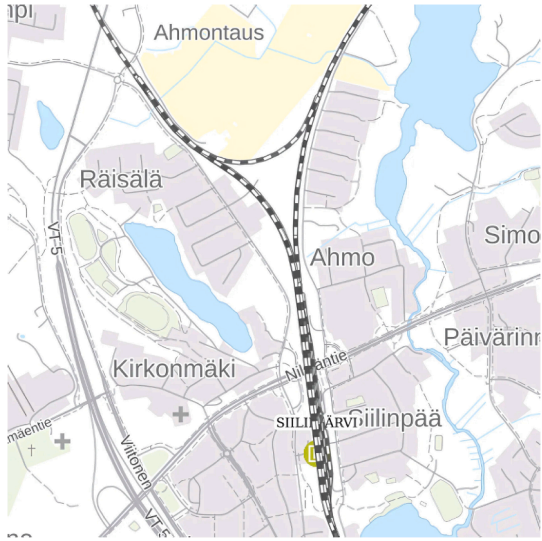

(a) Map

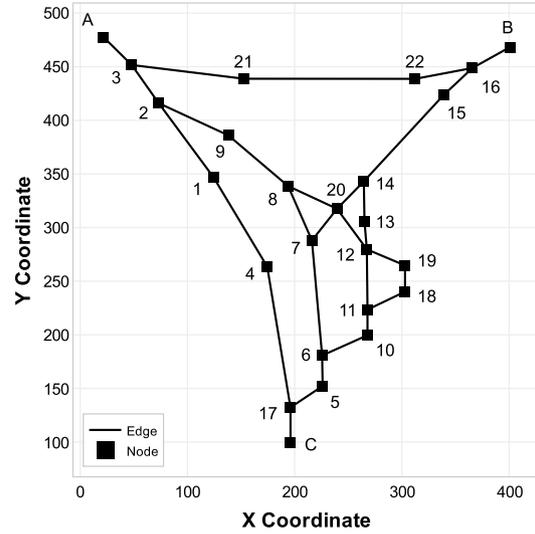

(b) Graph representation

**Fig. 3.** Representation of the Siilinjärvi railway station.

**Table 3**
Number of cost-efficient portfolios for different information sets and reliability requirements.

| Fortified switches | 1 | 3 | 5 | 10 | 15 | 20 |
| --- | --- | --- | --- | --- | --- | --- |
| **Feasible portfolios** | 23 | 1.7k | 35.4k | 1.7M | 4M | 4.2M |
| **Non-informative** $S^0$ | 5 | 41 | 163 | 972 | 1.7k | 1.9k |
| **Partial information** $S^1$ | 4 | 17 | 70 | 554 | 937 | 1067 |
| $S^1 \wedge \alpha = 0.95$ | 0 | 3 | 9 | 85 | 229 | 294 |
| $S^1 \wedge \alpha = 0.96$ | 0 | 0 | 2 | 24 | 119 | 178 |

with neural networks [50] or Dijkstra's algorithm [51], or deep neural networks [52], among others.

### 4.3. Preferences regarding connections

We consider two sets of preferences regarding the importance of operational connections. In the first, there is no preference information about the importance of connections. This situation is represented by the set $S^0 = \{w \in \mathbb{R}^3 \mid w_j \geq 0 \ \forall j, \sum_{j=1}^{3} w_j = 1\}$.

In the second, the connections are ranked based on the following yearly traffic volume: (A,B) at 3035 trains/year, (A,C) at 7547 trains/year, and (B,C) at 1373 trains/year. Specifically, the ratios between the connections' flow and the total flow are used to define the set of feasible weights. This gives the set $S^1 = \{w \in \mathbb{R}^3 \mid w_j \geq 0 \ \forall j, \sum_{j=1}^{3} w_j = 1, w_2 \geq 5.5 w_3, w_1 \geq 2.2 w_3\}$. For this situation, we also explore what fortification actions are needed to ensure that the reliability of each connection exceeds 95 % and 96 %. These reliability requirements are represented by $S^1 \wedge \alpha = 0.95$ and $S^1 \wedge \alpha = 0.96$, respectively.

### 4.4. Results

#### 4.4.1. Cost-efficient portfolios

We evaluated portfolios of up to 20 fortification actions. The corresponding cost-efficient portfolios, computed with Algorithm 1 for the two information sets and minimum reliability requirements, are in Table 3. For reference, the number of feasible portfolios is provided for different numbers of fortified switches.

When there are planned reliability requirements and preference information based on traffic volume, the number of cost-efficient portfo-lios becomes much smaller. Thus, there are fewer feasible portfolios of fortification actions, which helps to provide recommendations.

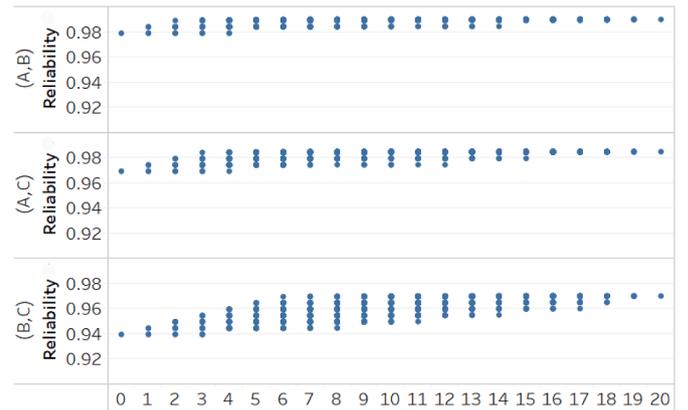

**Fig. 4.** Reliability of the connections for cost-efficient portfolios when there is no preference information about the importance of connections.

#### 4.4.2. Reliability of connections

The reliability of connections for cost-efficient portfolios based on the set $S^0$ is in Fig. 4. Each point corresponds to a different portfolio for a given connection. The marginal reliability improvement gained by fortifying an additional switch decreases with the number of actions that have already been implemented. This is because the most impactful actions are implemented first, and because the actions have equal costs. However, if the actions are not of equal cost, the marginal improvement may not decrease because there can be very impactful, high-cost actions that can only be implemented if there are sufficient resources to do so.

#### 4.4.3. Minimum performance for a single connection

Fig. 5 shows the reliability of connections for all cost-efficient portfolios containing up to five fortified switches for the cases of non-informative and partial preference information. In the case of partial information, there is relatively little improvement in the reliability of the connection (B, C) compared to the non-informative case, because this connection has less traffic than those that include the border node A.

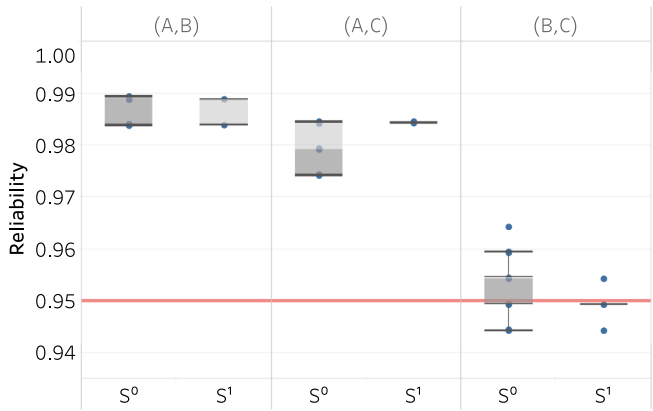

**Fig. 5.** Reliability of connections for cost-efficient portfolios for up to five fortified switches. Results are presented for the non-informative set ($S^0$) and the partial information set ($S^1$).

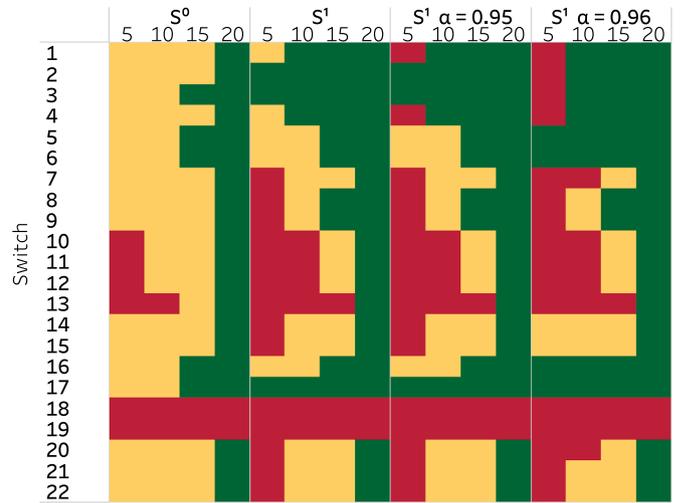

**Fig. 6.** Core index for the individual fortification actions given a different number of fortified switches and information sets. Red: core index 0, green: core index 1, and yellow: otherwise.

The non-informative set contains the weights (1,0,0), (0,1,0), and (0,0,1). These weights represent situations in which the focus is on a single connection; for example, for the weights (1,0,0), the reliability of the first connection is the only one that contributes to the utility function. Therefore, the maximum reliability that can be achieved for every single connection is given by cost-efficient portfolios for the non-informative set (see Lemma 1). For example, the maximum reliability of connection (B, C) that can be achieved by fortifying up to five switches is around 0.96. Moreover, if the minimum reliability requirement is set at 0.95, then some of the cost-efficient portfolios derived from the partial information set $S^1$ do not meet this requirement. If the reliability requirement is tightened to 0.96, this requirement can be met only by portfolios that are not cost-efficient in $S^1$.

The analysis of cost-efficient portfolios (see the last two rows of Table 3), computed with Algorithm 2 subject to reliability requirements at 0.95 and 0.96 for each connection, provides information about how many switches need to be fortified to fulfill reliability requirements. For example, at least three switches must be fortified to achieve a reliability level of 0.95 and five for 0.96.

*4.4.4. Selecting switches to fortify*

The composition of cost-efficient portfolios can be examined to derive recommendations for selecting fortification actions. If there is a single cost-efficient portfolio at a budget level (i.e., at a given number of fortified switches), the actions in this portfolio should be selected at this budget level. However, in general, there are many cost-efficient portfolios. In such cases, the *core index* of the actions can be computed to derive recommendations [53]. This index is the relative share of the non-dominated portfolios that contain the action. By definition, all cost-efficient portfolios with the same cost are non-dominated. We denote by $Q_{ND}(c)$ the set of non-dominated portfolios of cost $c$. Then, the core index of an action $l$ when there are $c$ resources is

$$CI(l,c) = \frac{|\{q \in Q_{ND}(c) \mid q_l = 1\}|}{|Q_{ND}(c)|}, \qquad (6)$$

where $|\cdot|$ is the cardinality of the set.

At a given budget level, a fortification action with a core index of 1 can be recommended because it is in every non-dominated portfolio. If the core index is 0, the action can be discarded. Eliciting additional preferences about the importance of connections between border nodes tends to reduce the number of cost-efficient portfolios and, therefore, change the core index of the actions, too.

Fig. 6 shows the core indices of the fortification actions for switches for the two different information sets with and without reliability requirements at different budget levels. If only a few switches can be fortified, there are fewer fortification actions with a core index of 1. This is because individual actions cannot usually improve all objectives, and when the budget is small, it is not possible to implement many actions, so that all objectives can be improved. At higher budget levels, which allow more switches to be fortified, some actions are in the portfolios that improve the reliability of all connections, and thus, their core index becomes 1.

If there are minimum performance requirements, there are far fewer cost-efficient portfolios that satisfy these constraints, too. Thus, the core indices tend to be 1 or 0. Some of the previously recommended actions will be discarded, such as fortifying switch 3 at a budget level of five switches. In contrast, some actions become relevant, such as fortifications of switches 16 and 17, even for small budgets. Thus, performance requirements should be introduced at the outset rather than afterward.

For comparison, centrality metrics such as *degree* or *betweenness* are widely used to quantify the relevance of nodes in a network. For example, Ghorbani-Renani et al. [54] examine five centrality metrics to select nodes to be protected in interdependent water, gas, and power networks. However, a concern about using centrality metrics to guide fortification actions is that these metrics do not necessarily reflect the performance of the networks. As we have shown, the prioritization of switches for fortification depends on the importance of operational connections and the minimum reliability requirements on the connections. For example, switch 10 is ranked highly by some centrality metrics (closeness 2nd, degree 3rd, and betweenness 6th), but it is only relevant for fortifying the network if many other switches can be fortified.

*4.4.5. Impact of estimated failure probability*

Although failure probabilities can be estimated from historical data or reliability models, this estimation task may be difficult due to reasons such as sparse or inaccurate data. Accordingly, it is instructive to explore how sensitive the recommendations are to these probability estimates. We thus assess how alternative assumptions about failure probabilities before ($p_k$) and after ($p'_k$) fortification affect recommendations for fortification actions. We consider[1] all combinations of the failure probabilities $p_k \in \{0.01, 0.02, 0.03, 0.04, 0.05\}$ and failure probabilities after fortification $p'_k \in \{p_k/2, p_k/3, p_k/4, p_k/5\}$. For each combination, all the nodes have the same parameters. The cost-efficient portfolios have the same composition for every combination.

---

[1] We do not present higher failure probabilities because the computation of the utility function is based on a cut set method, with lower precision for high probability values. Nevertheless, if necessary, the accuracy can be improved by adding more detail to this calculation.

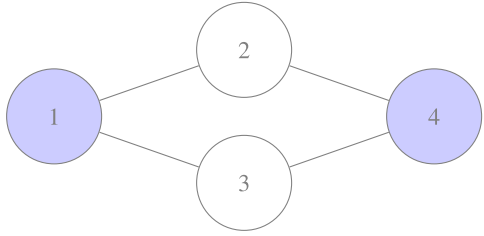

**Fig. 7.** Transportation network with border nodes: 1 and 4; and non-border nodes: 2 and 3.

A special case occurs for perfect fortification at $p'_k = 0$, which means that the possibility of disruption is completely eliminated at the node. In this case, for a given $p_k$, there are fewer cost-efficient portfolios than for $p'_k > 0$. If $p'_k = 0$, the cost-efficient portfolios are the same for $p_k \in \{0.01, 0.02, 0.03, 0.04, 0.05\}$.

To understand this, consider the network in Fig. 7, where switches 2 and 3 have a failure probability $p_2 = p_3 = 0.1$. The action portfolio is $q = (q_2, q_3)$. For the imperfect fortification ($p'_k > 0$), there are four cost-efficient portfolios: (0,0), (1,0), (0,1), (1,1). However, if the fortification is perfect, the fortification actions at switch 2 or switch 3 are sufficient to guarantee that there is a path between the border nodes. Then, $\mathbb{E}[u(x) \mid (1,0)] = \mathbb{E}[u(x) \mid (0,1)] = \mathbb{E}[u(x) \mid (1,1)] = 1$. As (1,1) costs more than (1,0) and (0,1), it is not cost-efficient.

The significance of this result is that recommendations about fortification actions can be provided even if the estimated failure probabilities are not fully accurate. Relatively accurate estimates would be needed to determine the network's performance or to compute cost-efficient portfolios that ensure minimum reliability requirements for individual connections.

## 5. Discussion

Our methodology assumes that individual node disruptions are independent of the disruptions to other nodes. In some situations, however, the disruption probability of a node may depend on other nodes. This would be the case, for example, if the failure of a node increases the load on another node, thereby increasing the disruption probability of the latter. Analyses of such interdependencies could be captured, for example, through Bayesian methods [see, e.g., [55,56]].

Challenges arising from the resulting increase in the number of required parameters can be alleviated by admitting partial information about disruption probabilities. This information could be obtained, for example, by asking experts to provide probability statements on verbal scales and mapping these statements to an interval of probabilities [see, e.g., [57]].

Using binary variables to model disruptions assumes that nodes are operational or disrupted. To model different levels of node performance more comprehensively, one can introduce multi-state reliability variables to capture different gradations of node performance and, for example, to model situations where the nodes can enable a reduced traffic volume [see, e.g., [58,59]]. A potential challenge with multi-state variables is that the number of network states grows quickly (e.g., if there are three rather than two states at each node in a network with $n$ nodes, the number of network states grows from $2^n$ to $3^n$).

Our approach can be linked to scenarios that characterize external operating conditions. For example, considering different weather scenarios could be instructive if the node disruption probabilities depend on the weather. This allows one to explore the robustness of cost-efficient portfolios by assessing which portfolios perform relatively well across all scenarios [39]. Moreover, if probabilities are associated with scenarios, one can gain insight into the expected performance of the network under those scenarios.

While our case study focuses on switches within a single station, the portfolio optimization approach can be applied to other network components, provided that their operational state can be represented with binary (operational/disrupted) variables. This approach can also be extended to study larger systems comprising several interconnected sub-networks or distinct layers of interconnected networks that provide a range of services. Yet, the resulting optimization problems are NP-hard even with completely specified weight information (solving them would imply solving the INLP formulation in Appendix A, which is NP-hard), and are therefore likely to become intractable for large networks. Advances in this area can be pursued by building hierarchical models (see, e.g., [60,61]) so that, for example, the regional or national railway network is subdivided into nodes representing railway stations, connected by railway tracks that correspond to edges. In this setting, cost-efficient portfolios at the station level can be aggregated to explore how resources should be allocated more generally to best contribute to the network's overall performance.

Our illustrative example considers a single fortification action per node. However, the same methodology can be used to analyze choices among multiple actions at a single node or even cross-cutting actions that affect the reliability of multiple nodes. From a modeling perspective, these extensions are relatively straightforward, as they can be accommodated through logical constraints. For example, suppose that a cross-cutting fortification action improves the reliability of multiple nodes. This is equivalent to our basic formulation, subject to the constraint that the corresponding improvement in reliability is attained at all the nodes affected by this action. An additional extension would be to explicitly consider the sequential dynamics of implementing fortification actions. That is, if there are constraints on how many fortification actions can be implemented per unit of time within the planning horizon, the question becomes *in which order* the actions should be implemented so that the performance of the network can be improved as quickly as possible. Importantly, one could also identify the optimal restoration order of the nodes, aiming to maximize the resilience of the networks after disruptions.

From the viewpoint of practice, DMs can use our methodology for several tasks in the maintenance and fortification of switches or other components of the railway system. For example, by identifying which components contribute most to the reliability of different connections, the methodology helps guide maintenance activities based on recent information about traffic volumes and disruption probabilities. Moreover, in planning monitoring activities, it is useful to know which components matter most for connections' reliability, so that more monitoring resources can be allocated to them. If the maintenance budget is to be reduced, the methodology can be used to identify which maintenance actions can be postponed with minimal compromise to the reliability of the connections.

## 6. Conclusion

In this paper, we have developed a multicriteria portfolio optimization approach to support the fortification of infrastructure networks, whose nodes may be disrupted by events such as natural hazards, technical failures, or intentional attacks. This optimization is based on a probabilistic risk assessment for quantifying disruption impacts and identifying which nodes of the network should be fortified, considering multiple objectives, resource constraints, and potential reliability requirements. Our approach explicitly accounts for the preferences for relevant networks' performance objectives, for example, the reliability of operational connections between designated pairs of nodes.

The proposed approach opens avenues for further research on strengthening critical infrastructure systems. A relatively straightforward extension is to consider infrastructure networks for the transportation of multiple commodities (e.g., multicommodity railway networks [62,63]), as such commodities can be prioritized in the same way as the connections between border nodes in this paper. A somewhat more demanding extension would be to consider multiple interdependent networks for energy, transportation, and communication, given that the dependencies between these networks would need to be modeled. A further extension would be to model situations in which the disruption probabilities are contingent not only on the selected fortification actions but also on the state of other nodes or even external conditions, as depicted by scenarios.

**Data availability**

Data will be made available on request.

**CRediT authorship contribution statement**

**Joaquín de la Barra:** Writing – review & editing, Writing – original draft, Visualization, Validation, Software, Methodology, Investigation, Formal analysis, Data curation, Conceptualization; **Ahti Salo:** Writing – review & editing, Writing – original draft, Validation, Supervision, Project administration, Methodology, Investigation, Funding acquisition, Formal analysis, Conceptualization; **Leevi Olander:** Writing – review & editing, Visualization, Validation, Software, Methodology, Investigation, Formal analysis, Data curation, Conceptualization; **Kash Barker:** Writing – review & editing, Validation, Supervision; **Jussi Kangaspunta:** Writing – original draft, Methodology, Investigation, Conceptualization.

**Declaration of interests**

The authors declare the following financial interests/personal relationships which may be considered as potential competing interests:

Leevi Olander reports financial support was provided by Finnish Transport Infrastructure Agency. Kash Barker reports financial support was provided by Fulbright Finland Foundation. If there are other authors, they declare that they have no known competing financial interests or personal relationships that could have appeared to influence the work reported in this paper.

**Acknowledgement**

This research has been partly supported by the Finnish Transport Infrastructure Agency (Väylävirasto), the Profi-8 project funded by the Academy of Finland, and a fellow award granted by the Fulbright Finland Foundation.

**Appendix A. Optimal fortification portfolio**

The solution to the nonlinear integer programming (INLP) problem in (A.1) gives the portfolios that maximize the expected utility when the preferences over the criteria, represented by the corresponding weights, are defined as exact values. In this formulation, one action per node is considered; therefore, the index $l$ is employed for both nodes and actions.

$$\underset{q}{\text{Maximize}} \quad \sum_{x \in \mathcal{X}} p(x \mid q) \sum_{j=1}^{m} w_j u_j(x) \quad \text{(A.1a)}$$

$$\text{Subject to:} \quad \sum_{l=1}^{h} q_l c_l \leq r \quad \text{(A.1b)}$$

$$Aq \leq b \quad \text{(A.1c)}$$

$$p_l = (1 - q_l) p_l^0 + q_l p_l', \quad l = 1, \ldots, h \quad \text{(A.1d)}$$

$$q_l \in \{0, 1\}, \quad l = 1, \ldots, h. \quad \text{(A.1e)}$$

The objective function (A.1a) represents the maximization of the expected utility, where the probabilities $p(x \mid q)$ of the network states $x$ depend on the fortification actions $q = (q_1, \ldots, q_h)$. Constraint (A.1b) ensures that the cost of implementing the selected portfolio does not exceed the budget $r$. Logical constraints are modeled by (A.1c) to represent mutually exclusive actions (e.g., $q_1 + q_2 \leq 1$) or actions which imply that some other action will also have to be implemented (e.g., $q_1 - q_2 \leq 0$); the matrix $A$ and vector $b$ contain the respective coefficients of such logical constraints. (A.1d) defines the dependency of the nodes' failure probabilities on the fortification actions. If the action is implemented at node $l$, then $q_l = 1$ and the failure probability is $p_l'$; otherwise, $q_l = 0$ and the failure probability is the initial failure probability $p_l^0$. Finally, (A.1e) represents fortification actions with binary decision variables.

**B. Upper bound of the portfolio performance**

The performance of Algorithm 1 depends on identifying portfolios that cannot be extended to cost-efficient ones. These portfolios should be discarded to avoid the complete enumeration of feasible portfolios. The earlier they are removed, the fewer portfolios that need to be evaluated. The method to compute the extended portfolios $q^a$ presented in the algorithm is simple and provides an exact upper bound (i.e., the maximum expected utility that can be achieved by a feasible portfolio constructed from the one under analysis) if the portfolios $q^a$ are feasible. On the other hand, if the portfolios are not feasible, the method can provide weak bounds.

An alternative approach to generate tighter bounds at a higher computational cost is presented. The solution of the optimization problem (B.1) provides an upper bound $\hat{u}_{w_e}$ of the expected utility at the extreme $w_e$, achievable by implementing a feasible portfolio by adding more fortification actions to portfolio $q$. Where $C$ and $D$ are the indices of the actions implemented and not implemented in the portfolio $q$, i.e., $C = \{l \mid q_l = 1\}$ and $D = \{l \mid q_l = 0\}$.

$$\begin{aligned}
&\underset{q}{\text{Maximize}} \quad u(q, w_e) \\
&\text{Subject to:} \quad \sum_{l=1}^{h} q_l c_l \leq r \\
&\quad q_l = 1, \quad \forall l \in C \\
&\quad q_l = 0, \quad \forall l \in D \\
&\quad q_l \in \{0, 1\}, \quad l = 1, \ldots, h.
\end{aligned} \quad \text{(B.1)}$$

Let $q^*$ denote a portfolio whose implementation provides the upper bounds of the expected utility at each extreme point, computed by solving (B.1). If any cost-efficient portfolio dominates portfolio $q^*$, the portfolio $q$ can be discarded because no cost-efficient portfolio can be generated from it. The main advantage of this approach is that the bounds are computed by considering the remaining available resources. One can also incorporate logical constraints into the optimization problem. Nevertheless, solving the problem (B.1) can lead to a different portfolio for each extreme point. If that is the case, it means that $q^*$ is not necessarily a feasible portfolio, and consequently, the upper bounds are not necessarily exact. Since $q^a$ implies implementing all the remaining fortification actions, the upper bounds provided by $q^*$ are always lower than or equal to (worst case) the ones provided by $q^a$.

The efficiency of the methods depends on several parameters, such as the utility function, available budget, and cost of the fortification actions, among others. For example, optimization problems with complex utility functions may be difficult to solve, making the computation of $\hat{u}_{w_e}$ impractical. In other cases, where the budget is small and just a few actions can be implemented, the bounds provided by $q_a$ can be very loose. One could also complement the methodologies by applying Algorithm 1 as presented and switch to the optimization model on the points where the portfolios $q^a$ are infeasible.